\documentclass[a4paper,10pt]{article}
\usepackage{amsfonts}
\usepackage{amsmath}
\usepackage{amsthm}
\usepackage{amsmath, graphics}
\newcommand{\mathsym}[1]{{}}

\begin{document}
\title{An application of Jacquet-Langlands correspondence to transfer operators for geodesic flows on Riemann surfaces}
\author{
Arash Momeni$^1$ and Alexei Venkov$^2$\\ 
$^1$ Department of Statistical Physics and Nonlinear Dynamics\\ Institute of Theoretical Physics, \\
Clausthal University of Technology, 38678 Clausthal-Zellerfeld, \\ Germany. E-Mail: arash.momeni@tu-clausthal.de \\
$^2$Institute for Mathematics, University of Aarhus, 8000 Aarhus C,\\ Denmark. E-Mail: venkov@imf.au.dk \\}
\maketitle

\begin{abstract}
In the paper as a new application of the Jacquet-Langlands correspondence we connect the transfer operators for different cofinite Fuchsian groups by comparing the corresponding Selberg zeta functions. 
\end{abstract}

\section{Transfer operator for cofinite Fuchsian groups and Selberg's zeta function}
\theoremstyle{plain}\newtheorem{th1}{Theorem}
\theoremstyle{plain}\newtheorem{th2}[th1]{Theorem}
\theoremstyle{plain}\newtheorem{th3}[th1]{Theorem}
\theoremstyle{plain}\newtheorem{th4}[th1]{Theorem}
\theoremstyle{plain}\newtheorem{th5}[th1]{Theorem}
\theoremstyle{plain}\newtheorem{th6}[th1]{Theorem}
\theoremstyle{plain}\newtheorem{th7}[th1]{Theorem}
\theoremstyle{plain}\newtheorem{th8}[th1]{Theorem}
\theoremstyle{plain}\newtheorem{th9}[th1]{Theorem}
\theoremstyle{plain}\newtheorem{th10}[th1]{Theorem}
\theoremstyle{remark}\newtheorem{rem}{Remark}
\ \ \ \ In this section we give a short summary of the generalization of Mayer's theory \cite{mayer91} by following Morita \cite{morita97}. In Mayer's theory a transfer operator $\mathcal L_s$ is introduced as a special case of Ruelle's operator for a dynamical system, for instance for the group $SL(2,\mathbb Z)$ this is the twice iterated Gauss map $T=T_G^2$ acting on the unit interval and the weight function $-\beta log\vert T'(z)\vert$. Then, as a result of the one-to-one correspondence between the closed geodesics on the modular surface 
$M=H\setminus SL(2,\mathbb Z)$ and the primitive periodic orbits of $T$ the Selberg zeta function can be written as a Fredholm determinant of the transfer operator \cite{mayer92}. 

 Morita in\cite{morita97} takes instead of $PSL(2,\mathbb R)$ and the Poincare half plane $H$ the isomorphic group $PSU(1,1)$ and as homogeneous space the unit disc $\mathbb D$ which is isomorphic to $H$. Let $\Gamma$ be a cofinite Fuchsian group. The canonical fundamental domain of $\Gamma$ is a polygon with a finite even number of sides $s_i$, each of which extends to a circular arc $C(s_i)$ perpendicular to the unit circle $S^1$, the boundary of $\mathbb D$. To every side $s_i$ a generator $g(s_i)$ is assigned which identifies the sides mutually. Now the action of the generators on the boundary points $S^1$ defines a Markov map \cite{bowen}
\begin{equation}
T_\Gamma x=g_ix\quad x\in S^1
\end{equation}
where one chooses for each $x$ the corresponding $g_i$ according to the location of $x$ between the footpoints of $C(s_i)$'s on $S^1$ as explained in \cite{morita97}. 

By letting the group $\Gamma$ act on sides $s_i$ and taking from them the arcs passing through vertexes of the fundamental polygon which cut the boundary, $S^1$ is partitioned into a set of intervals $\mathcal P'=\left\lbrace I_i \right\rbrace_i$. Like for the modular group and the Gauss map with the corresponding partition \cite{mayer book}, the map $T_\Gamma$ with respect to the partition $\mathcal P'$ satisfies the properties of orbit equivalence, piecewise monotonicity, the Markov property and transitivity. As in \cite{morita97} by using the intervals in $\mathcal P'$ it can be constructed a finite partition $\mathcal R=\left\lbrace A(1), A(2),..., A(q)\right\rbrace$  and its refinement $\mathcal P=\left\lbrace J\right\rbrace$ finite or infinite, such that the union of each of them differs from the other up to a set of zero Lebesgue measure. Furthermore to each $J\in \mathcal P$ there is a homomorphism $T_J$ from $J$ onto $A\in \mathcal R$. Finally a map $T$ is defined approximately everywhere on $X=\bigsqcup_{A\in \mathcal R}A$ such that $T\vert_{intJ}=T_J\vert_{intJ}$. As in \cite{morita97} we call the set of partitions $\mathcal R$ and  $\mathcal P$ together with the map $T$ a Markov system $\mathcal T_\Gamma=(\mathcal R, \mathcal P, T)$.

As in \cite{morita97} we set
\begin{equation}
G_J^s(w)=\vert T_J'(w)\vert^{-s}\quad w\in S^1, \quad s\in\mathbb C,
\end{equation}
then the transfer operator is defined as
\begin{equation}\label{transfer op.}
(L_{\mathcal T_\Gamma}(s)f)_i(z)=\sum_{J\in \mathcal P:\tau(J)=i}G_J^s(T_J^{-1}z)f_{\iota(J)}(T_J^{-1}z)
\end{equation}
on direct sum of $q= card \ \mathcal R$ Banach spaces chosen in such way that(\ref{transfer op.}) defines a nuclear operator where for $J\in \mathcal P$, $\iota (J)$ denotes the index with $int \ J\subset int \ A(\iota (J))$ and $\tau(J)$ denotes the index with $T_JJ=A(\tau(J))$. See the Convention 2.6 of \cite{morita97}. Like the case of $SL(2,\mathbb Z)$ regarding the correspondence between the closed geodesics on the fundamental surface of the group and the periodic orbit of $T$ we achieve a determinant expression for the corresponding Selberg zeta function. We summarize the final result as the following theorem of Morita \cite{morita97}. 
\begin{th1}
Let $\Gamma$ be a cofinite Fuchsian group. We denote $HC(\Gamma)$ to be the set of all primitive hyperbolic conjugacy classes in $\Gamma$. Then the Selberg zeta function $Z(s)$ has a determinant representation \begin{equation}
Z_{\Gamma}(s)\Xi(s)^2=Det(I-L_{\mathcal T_\Gamma}(s))
\end{equation}
where $\Xi(s)=\prod_{k=0}^\infty\prod_{c\in HC_1(\Gamma)}(1-exp(-(s+k)l(c)))$ for a certain finite subset $HC_1(\Gamma)$ of $HC(\Gamma)$. The function $\Xi(s)$ is meromorphic in $\mathbb C$ and, it is analytic in $\left\lbrace s\in\mathbb C:Re(s)>0\right\rbrace$ without zeros.
\end{th1}
We recall some notations and definitions we used in Theorm 1. According to Selberg for $Res>1$ the Selberg zeta function is given by absolutely convergent infinite product similar to Euler products of zeta function from number theory
\begin{equation}
Z_{\Gamma}(s)=\prod_{k=0}^\infty\prod_{\left\lbrace P\right\rbrace_\Gamma }(1-\mathcal N(P)^{-k-s})
\end{equation}
where $P$ runs through all primitive hyperbolic conjugacy classes $HC(\Gamma)$ with a norm $\mathcal N(P)>1$. $P$ primitive plays the role of prime number or prime ideal. Using his trace formula Selberg proved that \\ (1) $Z_{\Gamma}(s)$ has analytic(meromorphic) continuation to the whole $s\in \mathbb C$ \\ (2) $Z_{\Gamma}(s)$ satisfies the functional equation \begin{equation}
Z_{\Gamma}(1-s)=\Psi(s)Z_{\Gamma}(s)
\end{equation}
with known function $\Psi$ \\ (3) the nontrivial zeros of $Z_{\Gamma}(s)$ are related to eigenvalues and resonances of automorphic Laplacian for the group $A( \Gamma)$ \cite{selberg}, see also \cite{alexei}. 

There is another equivalent definition of the Selberg zeta function related to dynamical system of a geodesic flow on Riemann surface with constant negative curvature. This definition is important for theorem 1 above:
\begin{equation}
Z(s)=\prod_{k=0}^\infty\prod_{c\in HC(\Gamma)}(1-exp(-(s+k)l(c)))
\end{equation}
In formula above we use $c$ instead of $P$ like in \cite{morita97}. Recall that an element $c$ of $PSU(1,1)$ is hyperbolic if it has only two fixed point on the boundary $S^1$. $c$ is a primitive element of $\Gamma$ if it is not a power of other elements in $\Gamma$. Let $w_1,\ w_2\in \mathbb S^1$ are fixed points of $c$. This defines the geodesic line connecting $w_1,\ w_2$. The projection of this geodesic to the quotient space $\mathbb D\setminus \Gamma$ is closed geodesic and $l(c)$ is the hyperbolic length of corresponding prime closed geodesic. We have
\begin{equation}
\mathcal N(c)=e^{l(c)}
\end{equation}

\section{Jacquet-Langlands correspondence}
In \cite{hejhal} an explicit integral operator lift with Siegel theta function as the kernel, between Maass forms of unit group of an indefinite quaternion division algebra and congruence subgroups of the modular group is constructed. This is indeed an especial case of Jacquet-Langlands correspondence which Hejhal reproved by using classical argument \cite{hejhal}. To introduce this correspondence we follow \cite{hejhal}, \cite{bolte1} and \cite{bolte}. We start by recalling the unit group of quaternion algebra.
\subsection{Unit group of quaternion algebra}
In this part we follow \cite{miyake}.
We call a ring $B$ with unity an algebra of dimension $n$ over a field $F$, if the following three conditions are satisfied:

(1) $F\subset B$, and the unity of $F$ coincides with the unity of $B$;

(2) all elements of $F$ commute with all element of $B$;

(3) $B$ is a vector space over $F$ of dimension $n$. 
\\Let $B$ be an algebra over $F$. The center of $B$ is defined as the set of all commuting element of the algebra and is denoted by $Z(B)$. We call $B$ a central algebra when $Z(B)=F$. The algebra $B$ is called simple if it is simple as a ring, namely, if $B$ has no two-sided ideals except for $\left\lbrace 0\right\rbrace $ and $B$ itself. We call $B$ a division algebra if every nonzero element of $B$ is invertable. Now we can give the definition of quaternion algebra.\\A central simple algebra $B$ of dimension $4$ over a field $F$ is called a quaternion algebra over $F$. Furthermore, if $B$ is a division algebra, we call $B$ a division quaternion algebra. Let $B$ be a quaternion algebra over a field $F$. As a result of Wedderburn's theorem there are only two possibilities \cite{lewis}

(1) $B$ is a division quaternion algebra

(2) $B$ splits over $F$ that is, $B$ is  isomorphic to $M_2(F)$\\

Using the following results \cite{weil} one can define a norm for $B$. 1)If $F$ is algebraically closed, then $M_2(F)$ is a unique quaternion algebra over $F$ up to isomorphism. 2) Let $K$ be any extension over $F$. Then $B\otimes_FK$ is a quaternion algebra over $K$. We say that $B$ is ramified or splits over $K$ if $B\otimes_FK$ is a division quaternion algebra or is isomorphic to $M_2(K)$, respectively. Now for the quaternion algebra $B$ over $F$ let $\overline{F}$ be the algebraic closure of $F$. According to the results above $B\otimes_{F}\overline{F}$ is a quaternion algebra over the algebraically closed field $\overline{F}$ and therefore $B\otimes_{F}\overline{F}$ is isomorphic to $M_2(\overline{F})$. Now we can define (reduced) norm and (reduced) trace for elements of $B$ by
\begin{equation}
N_B(\beta)=det(\beta)\qquad tr_B(\beta)=tr(\beta)
\end{equation}
where $det(\beta)$ and $tr(\beta)$ are determinant and trace of $\beta$ as an element of $M_2(\overline{F})$. In \cite{weil} it is proved that both $N_B(\beta)$ and  $tr_B(\beta)$ belong to $F$.

For an algebra $B$ not necessarily of quaternion type over the field of rational numbers $F=\mathbb Q$ or its p-adic extensions $F=\mathbb Q_p$ the concept of an order is defined as a subset $R$ of $B$ if the following two conditions are satisfied:\\
(1) $R$ is a subring containing the unity of $B$\\
(2) $R$ is finitely generated over $\mathbb Z$ or $\mathbb Z_p$ and contains a basis of $B$ over $F$.\\An order of $B$ is called maximal if it is maximal with respect to inclusion. Conjugate orders are characterized by a number called discriminant. The discriminant $d(R)$ of an order $R$ of an algebra $B$ is defined to be the square of the product of primes $p$ where $B$ is ramified over $\mathbb Q$ that is $B\otimes_{\mathbb Q}\mathbb Q_p$ is division quaternion algebra. All orders with same discriminant are conjugate.\\Let $B$ be a quaternion algebra over $\mathbb Q$. We call $B$ indefinite or definite according as $B\otimes_{\mathbb Q}\mathbb R$ is isomorphic to $M_2(\mathbb R)$ or is a division quaternion algebra. Let $B$ be an indefinite quaternion algebra over $\mathbb Q$. We fix an isomorphism of $B\otimes_{\mathbb Q}\mathbb R$ onto $M_2(\mathbb R)$ and consider $B$ as an subalgebra of $M_2(\mathbb R)$ through this isomorphism. Then the norm $N_B(\beta)$ of an element $\beta$ of $B$ is nothing but the determinant of $\beta$ as a matrix, by definition. Let $R$ be an order of $B$. We put
\begin{equation}
R^1=\left\lbrace \beta\in R\ \vert\ N_B(\beta)=1\right\rbrace \subset GL_2^+(\mathbb R)
\end{equation}
and call it the unit group of norm $1$ of $R$. 
\begin{th3}
Let $B$ be an indefinite quaternion algebra over $\mathbb Q$, and $R$ be an order of $B$. Then $R^1$ is a Fuchsian group of the first kind. Moreover, if $B$ is a division quaternion algebra, then $R^1\backslash H$ is compact
\end{th3}
Because of this theorem from now on we restrict ourself to an order $\mathcal O$ of an indefinite quaternion algebra  over the field of rational numbers $\mathbb Q$ and its unit group $\mathcal O^1$. By definition we can fix an emmbeding
\begin{equation}
\sigma:\mathcal O^1\longrightarrow M_2(\mathbb R)
\end{equation}
Furthermore according to definition we have \begin{equation}\label{vector rep.}
\mathcal O\cong e_1\mathbb Z\oplus e_2\mathbb Z\oplus e_3\mathbb Z\oplus e_4\mathbb Z\cong \mathbb Z^4
\end{equation}
where $e_1$, $e_2$, $e_3$ and $e_4$ form a basis for $B$ over $\mathbb Q$. 

For this representation of the order $\mathcal O$ as a four dimensional vector space over $\mathbb Z$ the norm is defined as a quadratic form which realized as a $4\times 4$ symmetric matrix with respect to the basis of the indefinite quaternion algebra $B$:
\begin{equation}
\left[ S'\right]_{i,j}=e_i*e_j 
\end{equation}
where $*$ denotes the algebra product.

\subsection{Siegel theta function}
 In this part we introduce the Siegel theta function by following \cite{hejhal}, \cite{bolte1}, and \cite{bolte}. For a symmetric matrix $S\in GL_n(\mathbb R)$ the majorant $P$ is defined to be a positive definite matrix such that $PS^{-1}P=S$. \\ We fix $S$ to be
\begin{equation}
S=\left( \begin{array}{cccc}
0&0&0&1\\
0&0&-1&0\\
0&-1&0&0\\
1&0&0&0\\
\end{array}\right) 
\end{equation}
It is easy to see that identity matrix is one of the majorants of $S$. \\ 
For $L_1,\ L_2\in SL_2(\mathbb R)$ and $m_1,\ m_2 \in M_2(\mathbb R)$ such that 
\begin{equation}
m_1=L_1m_2L_2^{-1}
\end{equation}
we define $A(L_1,L_2)\in M_4(\mathbb R)$ by
\begin{equation}
\left( \begin{array}{cccc}
\alpha_1\\
\beta_1\\
\gamma_1\\
\delta_1\\
\end{array}\right)=A(L_1,L_2)\left( \begin{array}{cccc}
\alpha_2\\
\beta_2\\
\gamma_2\\
\delta_2\\
\end{array}\right)
\end{equation}
where
\begin{equation}
m_1= \left( \begin{array}{cccc}
\alpha_1&\beta_1\\
\gamma_1&\delta_1\\
\end{array}\right), \ m_2=\left( \begin{array}{cccc}
\alpha_2&\beta_2\\
\gamma_2&\delta_2\\
\end{array}\right)
\end{equation}

As in \cite{hejhal} one can see that $A(L_1,L_2)^tA(L_1,L_2)$ is a majorant of $S$. For $w$ and $z$ belonging to $H$ where $z=x+iy$  we define
\begin{equation}
M_z:=\left( \begin{array}{cccc}
y^{\frac{1}{2}}&xy^{-\frac{1}{2}}\\
0&y^{-\frac{1}{2}}\\
\end{array}\right)\ \text{and}\ P_zw:=A(M_z^{-1},M_w^{-1})^tA(M_z^{-1},M_w^{-1})
\end{equation}
Let $\mathcal O$ be an order in an indefinite quaternion algebra over $\mathbb Q$. Then let $S'$ be the matrix of the norm of $\mathcal O$ with respect to a fixed basis of the algebra. For $q\in \mathcal O$ let $k_q\in Z^4$ be the corresponding element in $\mathbb Z^4$ in this basis (see formula (\ref{vector rep.}) from above subsection). 

We fix an embedding $\sigma:\mathcal O\longrightarrow M_2(\mathbb R)$. Since $\sigma$ is linear, we have a unique $B\in GL_2(\mathbb R)$ which for every $q$ satisfies 
\begin{equation}
\left( \begin{array}{cccc}
\alpha\\
\beta\\
\gamma\\
\delta\\
\end{array}\right)=Bk_q,\text{whenever}\ \sigma_q:=\sigma(q)=\left( \begin{array}{cccc}
\alpha&\beta\\
\gamma&\delta\\
\end{array}\right)
\end{equation}
Since
\begin{equation}
k_q^tS'k_q=2n(q)=2\ det(\sigma_q)=2(\alpha\delta-\beta\gamma)=(Bk_q)^tS(Bk_q)
\end{equation}
we conclude that $S'=B^tSB$. For this fixed embedding of $\mathcal O$, we define majorants $P_{zw}'$ of $S'$ by $P_{zw}':=B^tP_{zw}B$.\
Now fix $z_0\in H$ and let $\tau=u+iv,\ z=x+iy\in H$. With $R:=uS'+ivP_{zz_0}'$, we define a Siegel theta function $\theta(z;\tau)$ by
\begin{equation}
\theta(z;\tau):=Im(\tau) \sum_{k\in \mathbb Z^4}e^{\pi i k^tRk}=Im(\tau)\sum_{q\in \mathcal O}e^{\pi ik_q^tRk_q}
\end{equation}
The Siegel theta function has the following transformation property which is crucial for the application in the next subsection. We give it as a theorem whose proof one can find in \cite{bolte1} 
\begin{th4}
Let $\mathcal O$ be an order in an indefinite quaternion algebra over $\mathbb Q$, with (reduced) discriminant $d$. Then\\
$ (1)\ \theta(\sigma_qz;\tau)=\theta(z;\tau),\ \forall q\in \mathcal O^1$ \\
$(2)\ \theta(z;g\tau)=\theta(z;\tau),\ \forall g\in \Gamma_0(d)$
\end{th4}
\subsection{Theta-lifts}
In this subsection we introduce two integral transformations providing us a lift between Maass forms of unit group and congruence subgroup.\
As before, let $\mathcal O$ be an order with discriminant $d(\mathcal O)$ in an indefinite quaternion division algebra over $\mathbb Q$ and $\mathcal O^1$ denotes the corresponding unit group. In \cite{bolte1} by choosing a suitable referencee point $z_0$ appearing in the definition of Siegel theta function we have the following result which we give as a theorem.
\begin{th5}
For a suitably chosen reference point in the Siegel theta function the following maps 
\begin{equation}
\Theta:L_0^2(X_{\mathcal O})\longrightarrow \mathcal C_{d(\mathcal O)}\
\text{and} \ \overset{\sim}{\Theta}:\mathcal C_{d(\mathcal O)}\longrightarrow L_0^2(X_{\mathcal O})
\end{equation}
define bounded linear operators preserving Laplace eigenvalues given by
\begin{equation}
\Theta\varphi(\tau):=\int_{\mathcal F_{\mathcal O^1}}\theta(z;\tau)\varphi(z)d\mu(z)
\end{equation}
and
\begin{equation}
\overset{\sim}{\Theta}g(z):=\int_{F_d}\overline{\theta(z;\tau)}g(\tau)d\mu(\tau)
\end{equation}
where $L_0^2(X_{\mathcal O})$ and $\mathcal C_{d(\mathcal O)}$ denote the space of cusp forms of unit quaternion group $\mathcal O^1$ and congruence subgroup $\Gamma_0(d(\mathcal O))$ respectively. $\mathcal F_{\mathcal O^1}$ and $F_d$ are the corresponding fundamental domains.
\end{th5}
In \cite{bolte1} it is proved that for a basis $\left\lbrace  \varphi_k\ k\in\ \mathcal N_0\right\rbrace$ of eigenfunctions of Laplacian on $X_{\mathcal O}$, $\Theta(\varphi_k)$ is not identically zero for all $k$. Thus we have the following theorem whose proof one can see in \cite{bolte1}
\begin{th6}
All eigenvalues of the Laplacian on $L_0^2(X_{\mathcal O})$ also occur as eigenvalues of the Laplacian on $L^2(X_d)$
\end{th6}
In \cite{bolte} equality of the right hand sides of Selberg's trace formula for $\mathcal O^1$ and Selberg's trace formula for new forms of congruence subgroups which are defined in the next section, is shown. In this way one gets the following identity
\begin{equation}
\sum_{\varphi_k\in L^2(X_{\mathcal O})}h(r_k)=\sum_{g_k\in\mathbb C\oplus \mathcal C_d^{new}}h(r_k)
\end{equation}
This together with the last theorem leads to 
\begin{th7}
The eigenvalues of the Laplacian, including multiplicities, on $X_{\mathcal O}$ coincides with Laplace spectrum of Maass-new forms for the congruence subgroup $\Gamma_0(d)$, where the $d$ is the discriminant of the order $\mathcal O$.
\end{th7}
\section{Selberg trace formula}
In this section we introduce shortly the Selberg trace formula for the general case of a cofinite group together with the special cases of congruence subgroups and unit group of quaternion algebra. Then we define the Selberg trace formula for new forms.

The trace formula is a general identity connecting geometrical and spectral terms
\begin{equation}
\sum\left\lbrace \text{spectral terms}\right\rbrace =\sum\left\lbrace \text{geometric terms}\right\rbrace 
\end{equation}
In the Selberg trace formula the spectral terms come from the discrete and the continuous spectrum of automorphic hyperbolic Laplacian $A(\Gamma)$ for a cofinite group $\Gamma$ and the geometrical terms are integral operators depending on the conjugacy classes of $\Gamma$. As in \cite{alexei} one can calculate the integrals explicitly and achieve the final form of the Selberg trace formula. We give it in the form of a theorem as in \cite{koyama}
\begin{th2}
Let $h(r^2+\frac{1}{4})$ be a function of a complex variable $r$ which satisfies the following assertions:

(1) As a function of $r$, $h(r^2+\frac{1}{4})$ is holomorphic in the strip \\$\left\lbrace r\in \mathbb C:\vert Im(r)\vert<\frac{1}{2}+\varepsilon\right\rbrace$ for some $\varepsilon>0$. 

(2) In that strip, $h(r^2+\frac{1}{4})=O((1+\vert r^2\vert)^{-1-\varepsilon})$ and all of the series and integrals appearing below converge absolutely.

Then the following identity hold
\begin{equation}\label{trace formula}
\sum_{k=0}^{\infty}h(\lambda_k)+C(s)=I(s)+H(s)+E(s)+P(s)
\end{equation}
where $\left\lbrace \lambda_n \ \vert\ 0=\lambda_0<\lambda_1\leq \lambda_2\leq \ldots\right\rbrace$ are the discrete eigenvalues of $A(\Gamma)$. $C(s)$ corresponds to the continuous part of the spectrum given by
\begin{equation}
C(s)=-\frac{1}{4\pi}\int_{-\infty}^{\infty}\dfrac{\varphi'}{\varphi}(\frac{1}{2}+ir)h(r^2+\frac{1}{4})dr
-\dfrac{K_0}{4}h(\dfrac{1}{4})
\end{equation}
where $\varphi$ denotes determinant of the scattering matrix $\Phi(s)$ and $K_0:=\underset{s\rightarrow \frac{1}{2}}{lim}tr(\Phi(s))$. On the right hand side of (\ref{trace formula}) we have $I(s)$ which corresponds to the identity element of the group and given by
\begin{equation}
I(s)=\dfrac{vol(\Gamma\setminus H)}{4\pi}\int_{-\infty}^\infty r \tanh(\pi r)\ h(r^2+\frac{1}{4})dr
\end{equation}
The term $H(s)$ denotes the contribution of hyperbolic conjugacy classes and is given by
\begin{equation}
\sum_{\left\lbrace P\right\rbrace _\Gamma}\sum_{m=1}^\infty \dfrac{logN(P)}{N(P)^{\frac{m}{2}}-N(P)^{-\frac{m}{2}}}g(m \ logN(P))
\end{equation}
where $\left\lbrace P\right\rbrace _\Gamma$ denotes the primitive hyperbolic conjugacy class and the function $g$ appears in the process of the Selberg transformation:
\begin{equation}
g(u)=\frac{1}{2\pi}\int_{-\infty}^{\infty} e^{-iru}h(r^2+\frac{1}{4})dr
\end{equation}
The next term $E(s)$ refers to the contribution of elliptic elements given by a summation over primitive elliptic conjugacy classes $\left\lbrace R\right\rbrace _\Gamma$ with order $\nu$ 
\begin{equation}
E(s)=\frac{1}{2}\sum_{\left\lbrace R\right\rbrace _\Gamma}\sum_{m=1}^{\nu-1}\dfrac{1}{\nu \sin\pi m/\nu}\int_{-\infty}^\infty \dfrac{exp(-2\pi rm/\nu)}{1+exp(-2\pi r)}h(r^2+\frac{1}{4})dr
\end{equation}
Finally the last term comes from the parabolic conjugacy classes given by
\begin{equation}
P(s)=-K\ g(0)\ log2
+\dfrac{K}{4}h(\frac{1}{4})
-\dfrac{K}{2\pi}\int_{-\infty}^\infty\psi(1+ir)h(r^2+\frac{1}{4})dr
\end{equation}
where $K$ is the number of cusps and $\psi$ is the di-gamma function.
\end{th2}

In the case of cocompact groups including the unit group of quaternion algebras there is no continuous spectrum and no parabolic element, thus the trace formula (\ref{trace formula}) reduces to 
\begin{equation}
\sum_{k=0}^{\infty}h(\lambda_k)=I(s)+H(s)+E(s)
\end{equation}
 
Consider a congruence subgroup $\Gamma_0(n)\subset SL(2,\mathbb Z)$. As in \cite{bolte} the space $\mathcal C_n(\lambda)$ of Maass forms with eigenvalue $\lambda$, can be decomposed into two subspaces of new and old forms $\mathcal C_n(\lambda)=\mathcal C_n^{old}(\lambda)\oplus \mathcal C_n^{new}(\lambda)$. The space $\mathcal C_n^{old}(\lambda)$ is the linear span of all forms with the same eigenvalue $\lambda$ coming from all overgroups $\Gamma_0(m)\supset\Gamma_0(n)$ with $m\vert n$ and $ \mathcal C_n^{new}(\lambda)$ is defined to be the orthogonal complement of $\mathcal C_n^{old}(\lambda)$. Let us denote the dimension of $\mathcal C_n(\lambda)$ and $\mathcal C_n^{new}(\lambda)$ by $\delta(n,\lambda)$ and $\delta^{new}(n,\lambda)$ respectively. Then the following identity holds \cite{atkin}
\begin{equation}\label{dim-formula}
\delta^{new}(n,\lambda)=\sum_{m\vert n}\beta(\frac{n}{m})\delta(m,\lambda)
\end{equation}
with\begin{equation}
\beta(a)=\sum_{l\vert a}\mu(l)\mu(\frac{a}{l}),
\end{equation}
where $\mu(a)$ is the Moebius function. Identity (\ref{dim-formula}) leads to the following formula \cite{bolte}
\begin{equation}
\sum_{u_k\in \mathcal C_n^{new}} h(\lambda_k)=\sum_{m\vert n}\beta(\frac{n}{m})\sum_{u_k\in \mathcal C_m}h(\lambda_k)
\end{equation}
This suggest us to take the sum  
\begin{equation}
\sum_{u_k\in \mathcal C_n^{new}} h(\lambda_k)+\sum_{m\vert n}\beta(\frac{n}{m})h(\lambda_0)=\sum_{m\vert n}\beta(\frac{n}{m})\sum_{u_k\in \mathbb C\oplus\mathcal C_m}h(\lambda_k)
\end{equation}
for defining the left hand side of Selberg trace formula for new forms. Note, that under certain condition for $d$, as has been explained in \cite{bolte}, the sum $\sum_{m\vert n}\beta(\frac{n}{m})$ is equal to one and therefore  
\begin{equation}
\sum_{u_k\in \mathbb C\oplus\mathcal C_n^{new}} h(\lambda_k)=\sum_{u_k\in \mathcal C_n^{new}} h(\lambda_k)+h(\lambda_0)=\sum_{m\vert n}\beta(\frac{n}{m})\sum_{u_k\in \mathbb C\oplus\mathcal C_m}h(\lambda_k)
\end{equation}
Thus we achieve the following trace formula for new forms
\begin{equation}
\sum_{u_k\in \mathcal C_n^{new}} h(\lambda_k)+\sum_{m\vert n}\beta(\frac{n}{m})h(\lambda_0)=I^{new}(s)+H^{new}(s)+E^{new}(s)+CP^{new}(s)
\end{equation}The terms  $I^{new}(s)$, $H^{new}(s)$, $E^{new}(s)$ and $CP^{new}(s)$ are given by
\begin{equation}
I^{new}(s)=\sum_{m\vert n}\beta(\frac{n}{m})I_m(s)
\end{equation}
\begin{equation}\label{f35}
H^{new}(s)=\sum_{m\vert n}\beta(\frac{n}{m})H_m(s)
\end{equation}
\begin{equation}
E^{new}(s)=\sum_{m\vert n}\beta(\frac{n}{m})E_m(s)
\end{equation}
\begin{equation}
CP^{new}(s)=\sum_{m\vert n}\beta(\frac{n}{m})(P_m(s)-C_m(s))
\end{equation}
where $I_m(s)$, $H_m(s)$, $E_m(s)$ and $P_m(s)$ denotes the identity, hyperbolic, elliptic and parabolic conjugacy classes for the corresponding group $\Gamma_0(m)$. The term $C_m(s)$ refers to the contribution of the continuous spectrum for the corresponding group $\Gamma_0(m)$.
Now we give the following theorem which is proved in \cite{bolte}.
\begin{th8}
Let $\mathcal O$ be an order with discriminant $d$ in an indefinite quaternion division algebra over the field of rationals. Then the right hand side of the Selberg-new trace formula for congruence subgroup coincide with the right hand side of the Selberg trace formula for the unit group of $\mathcal O$
\begin{equation}\label{fff}
H_{\mathcal O^1}(s)=H^{new}(s)=\sum_{m\vert n}\beta(\frac{n}{m})H_m(s)
\end{equation}
where $H_{\mathcal O^1}(s)$ denotes the contribution of hyperbolic elements in the right hand side of Selberg trace formula for the unit group of quaternion.
\end{th8}
\section{Selberg zeta function and its determinant expression}
In this section we define Selberg zeta function for new forms and using the result of previous section we show that it is equal to Selberg zeta function of unit group of quaternion $\mathcal O^1$.\
As in \cite{koyama} we define the Selberg zeta function as the logaritmic derivitive of the hyperbolic contribution of the right hand side of Selberg trace formula, that is
\begin{equation}\label{Selberg zeta}
\dfrac{d}{ds}H(s)=\dfrac{d}{ds}\dfrac{1}{2s-1}\dfrac{d}{ds}log\ Z(s)
\end{equation}
Now by formula (\ref{f35}) we can define the Selberg-new zeta function as 
\begin{equation}
\dfrac{d}{ds}\sum_{m\vert n}\beta(\frac{n}{m})H_m(s)=\dfrac{d}{ds}H^{new}(s)=\dfrac{d}{ds}\dfrac{1}{2s-1}\dfrac{d}{ds}log\ Z^{new}(s)
\end{equation}
By simple calculation the later leads to the following formula
\begin{equation}
Z^{new}(s)=\Pi_{m\vert n}Z_m^{\beta(\frac{n}{m})}(s)
\end{equation}
where $Z_m(s)$ denotes the Selberg zeta function for congruence subgroup $\Gamma_0(m)$.
On the other hand theorem ($8$) together with (\ref{Selberg zeta}) leads to 
\begin{equation}
Z_{\mathcal O^1}(s)=Z^{new}(s)
\end{equation}
Let summerise our first result as a theorem
\begin{th9}
The following formula holds between different Selberg zeta functions
\begin{equation}
Z_{\mathcal O^1}(s)=Z^{new}(s)=\Pi_{m\vert n}Z_m^{\beta(\frac{n}{m})}(s)
\end{equation}
\end{th9}
Finally we can connect transfer operators for different groups through their Fredholm determinants and this is our second result in this paper. The detailed proof we publish later.
\begin{th10}
The following identity holds
\begin{equation}\label{final}
Det(I-L_{\mathcal O^1}(s))=\dfrac{\Xi_{\mathcal O^1}^2(s)}{\Pi_{m\vert n}\Xi_{m}^{2\beta(\frac{n}{m})}(s)}\Pi_{m\vert n}Det(I-L_{\Gamma_0(m)}(s))^{\beta(\frac{n}{m})}
\end{equation}

\end{th10} 
\begin{rem}
 The right hand side of (\ref{final}) can be represented in a different form. That is
\begin{equation}
\Xi_{\mathcal O^1}^2(s) \Pi_{m\vert n}Det(I-\mathcal L_s^{\Gamma_0(m)})^{\beta(\frac{n}{m})}
\end{equation}
Here the operator $\mathcal L_s^{\Gamma_0(m)}$ was introduced and studied by Mayer (see for example \cite{chang}). The operator acts on a certain Banach space of vector valued functions $f:\Omega\rightarrow \mathbb C^{\mu}$, where $\Omega$ is a domain in $\mathbb C$ and $\mu$ is the index of $\Gamma_0(m)$ in $SL(2,\mathbb Z)$. 
 
\end{rem}

\textbf{Acknowledgement.} In this work the first author was supported by German Academic Exchange Service (DAAD) and the second author was supported by Alexander von Humboldt Foundation. \\ Both authors would like to thank Dieter Mayer for important remarks(including remark.1), his hospitality and providing excelent work's condition in his institute.

\end{document}